\documentclass[12pt]{article}

\usepackage{xspace}
\usepackage{url}
\usepackage{mathtools}
\usepackage{amssymb}
\usepackage{amsthm}
\usepackage{empheq}
\usepackage{latexsym}
\usepackage{enumitem}
\usepackage{eurosym}
\usepackage{dsfont}
\usepackage{appendix}
\usepackage{color} 
\usepackage[unicode]{hyperref}
\usepackage{frcursive}
\usepackage[utf8]{inputenc}
\usepackage[T1]{fontenc}
\usepackage{geometry}
\usepackage{multirow}
\usepackage{todonotes}
\usepackage{lmodern}
\usepackage{anyfontsize}  
\usepackage{stmaryrd}
\usepackage{cleveref}
\usepackage[english]{babel}
\usepackage[english=british]{csquotes}
\usepackage{mathtools}
\usepackage{accents}
\usepackage{color}
\usepackage{comment}
\usepackage{bm}

\definecolor{red}{rgb}{0.7,0.15,0.15}
\definecolor{green}{rgb}{0,0.5,0}
\definecolor{blue}{rgb}{0,0,0.7}
\hypersetup{colorlinks, linkcolor={red}, citecolor={green}, urlcolor={blue}}

\allowdisplaybreaks

\usepackage{natbib}
\setcitestyle{numbers,open={[},close={]}}

\definecolor{red}{rgb}{0.7,0.15,0.15}
\definecolor{green}{rgb}{0,0.5,0}
\definecolor{blue}{rgb}{0,0,0.7}
\hypersetup{colorlinks, linkcolor={red}, citecolor={green}, urlcolor={blue}}
			
\makeatletter \@addtoreset{equation}{section}

\newtheorem{theorem}{Theorem}[section]
\newtheorem{assumption}[theorem]{Assumption}

\newtheorem{example}[theorem]{Example}

\newtheorem{lemma}[theorem]{Lemma}
\newtheorem{proposition}[theorem]{Proposition}

\newtheorem{definition}[theorem]{Definition}
\newtheorem{remark}[theorem]{Remark}

\DeclareUnicodeCharacter{014D}{\=o}
\newcommand{\ba}{\begin{array}}
\newcommand{\ea}{\end{array}}
\newcommand{\be}{\begin{equation}}
\newcommand{\ee}{\end{equation}}
\newcommand{\bea}{\begin{eqnarray}}
\newcommand{\eea}{\end{eqnarray}}
\newcommand{\beaa}{\begin{eqnarray*}}
\newcommand{\eeaa}{\end{eqnarray*}}

\def\dbA{\mathbb{A}}

\def\dbE{\mathbb{E}}
\def\dbF{\mathbb{F}}

\def\dbI{\mathbb{I}}

\def\dbL{\mathbb{L}}

\def\dbN{\mathbb{N}}
\def\dbP{\mathbb{P}}
\def\dbR{\mathbb{R}}

\def\a{\alpha}

\def\g{\gamma}
\def\d{\delta}

\def\l{\lambda}
\def\m{\mu}
\def\n{\nu}
\def\si{\sigma}
\def\t{\tau}
\def\f{\varphi}
\def\th{\theta}

\def\D{\Delta}
\def\Th{\Theta}

\def\O{\Omega}
\def\cA{{\cal A}}
\def\cB{{\cal B}}

\def\cD{{\cal D}}

\def\cF{{\cal F}}

\def\cH{{\cal H}}

\def\cJ{{\cal J}}

\def\cL{{\cal L}}
\def\cM{{\cal M}}

\def\cP{{\cal P}}

\def\cW{{\cal W}}

\def\q{\quad}
\def\pa{\partial}

\def\qed{ \hfill \vrule width.25cm height.25cm depth0cm\smallskip}

\def\bx{{\mathbf{x}}}

\def\1{\mathbf{1}}

\def\bS{{\mathbf{S}}}

\begin{document}

\title{\bf{It\^o-Wentzell formulas for \\ semimartingale conditional laws \\ with applications to mean-field control}}\author{Mehdi Talbi\footnote{Laboratoire de Probabilités, Statistiques et Modélisation, Université Paris-Cité, France, talbi@lpsm.paris} \quad Nizar Touzi\footnote{Tandon School of Engineering, New York University, USA, nizar.touzi@nyu.edu. This author is supported in part by NSF grant \#DMS-2508581. }} 
\date{\today}

\maketitle

\begin{abstract}
The present paper is an extension of \cite{fadle2024ito}. Following the same methodology, merely based on Taylor expansions, we establish the Itô and Itô-Wentzell formulae for flows of conditional distributions of general semimartingales, thus allowing for discontinuous semimartingales with possibly discontinuous flows of conditional marginals. We apply these results to derive the dynamic programming equations corresponding to mean field control problems with Poisson type common noise and mean field stopping problems with common noise.
\end{abstract}

\textbf{MSC2020.} 60G40, 49N80, 35Q89, 60H30. \\

\textbf{Keywords.} Itô formula, Itô-Wentzell formula Wasserstein space, mean field optimal control, mean field optimal stopping.

\section{Introduction}

This paper is the continuation of \cite{fadle2024ito} who derive an Itô and an Itô-Wentzell formulas for flows of conditional laws of  continuous semimartingales. Our main objective here is to extend the same methodology to the general context of càdlàg semimartingales with possibly discontinuous flow of conditional marginals. One of the highlights of this work consists in its methodology, which is essentially based on the use of the linear functional derivative and Taylor expansions, in the same spirit as the proof of the traditional finite dimensional Itô formula. 

Due to the rising interest for mean field games and related problems, independently introduced by \citeauthor*{lasry2007mean} \cite{lasry2007mean} and \citeauthor*{huang2006large} \cite{huang2006large}, chain rules for flows of marginals have received a strong attention over the past few years, and have been established through diverse methods. The most popular approach consists in establishing the formula for the $N$-particles approximation of the mean field system, and letting $N$ go to infinity to derive the mean field Itô formula; see e.g.\ \citeauthor*{chassagneux2022probabilistic} \cite{chassagneux2022probabilistic}, \citeauthor*{buckdahn2014mean} \cite{buckdahn2014mean} or \citeauthor*{carmona2018probabilisticI} \cite{carmona2018probabilisticI, carmona2018probabilisticII}. Another method relies on the density of cylindrical functions in the space of smooth functions, see \citeauthor*{guo2023ito} \cite{guo2023ito} and \citeauthor*{guo2024ito} \cite{guo2024ito}. In \cite{talbi2023dynamic}, \citeauthor*{talbi2023dynamic} exploit the properties of the linear functional derivative and the standard Itô formula to derive a general Itô formula for deterministic flows of càdlàg semimartingales. 

In the context of random fields, we may mention the contribution of \citeauthor*{dos2023ito} \cite{dos2023ito}, who proved an Itô-Wentzell formula for continuous conditional flows of semimartingale measures by using the finite population approximation. The recent paper by \citeauthor*{jisheng2025wentzell} \cite{jisheng2025wentzell} proves this formula for general càdlàg semimartingales, relying on the approximation of the characteristics of the random field by cylindrical functions.

In the present paper, we show that the methodology of \cite{talbi2023dynamic} and \cite{fadle2024ito} can be used in the context of a càdlàg semimartingale $X$. More precisely, given a deterministic mapping $t \mapsto u(m_t)$ or a random field $t \mapsto U_t(m_t)$, $m_t := \dbP \circ X_t^{-1}$, and using a partition $\pi := \{t_0 := 0, t_1, \dots, T_N := T\}$ of the time interval $[0,T]$, we rewrite the decomposition
$$ u(m_T) = \sum_{n=1}^N u(m_{t_n}) - u(m_{t_{n-1}}) $$
with the linear functional derivative $\d_m u$ to obtain Taylor-like expansions on each subinterval $(t_{n-1}, t_n]$, which enables us to derive the formula by letting the step of the partition tend to $0$. Here, the main difficulty compared to the continuous case is that the remainder in Taylor expansions does not necessarily vanish as the step of the partition goes to $0$, due to the presence of jumps. We overcome this difficulty by adapting the proof of the finite dimensional generalized Itô's formula to our context. 

The paper is structured as follows. In Section \ref{sec:Ito}, we establish the general It\^o formula for conditional flows of semimartingale measures. In Section \ref{sec:Ito-Wentzell}, we extend this formula to the case where the function of the flow of marginals is a random field, which corresponds to the Itô-Wentzell formula. Section \ref{sec:examples} illustrates the usefulness of our results with two examples from mean field stochastic control problems, namely, mean field control with Poisson type common noise and mean field optimal stopping with common noise. Finally, Appendix \ref{sec:appendix} contains the proofs of important technical lemmas. 

\medskip
\noindent \textbf{Notations.} We denote by $\cP(\O,\cF)$ the set of probability measures on a measurable  space $(\O,\cF)$, and $\cP_2(\O,\cF)$ the subset of square integrable probability measures in $\cP(\O,\cF)$, equipped with the $2$-Wasserstein distance $\cW_2$. When $(\O,\cF) = (\dbR^d,\cB(\dbR^d))$, we simply denote them as $\cP(\dbR^d)$ and $\cP_2(\dbR^d)$. $\cD^d$ denotes the space of càdlàg paths from $[0,T]$ to $\dbR^d$. For a random variable $Z$ and a probability $\dbP$, we denote by $\dbP_Z:=\dbP\circ Z^{-1}$ the law of $Z$ under $\dbP$. 
 For vectors $x, y\in \dbR^n$ and matrices $A, B\in \dbR^{n\times m}$, denote $ x\!\cdot\! y:=\sum_{i=1}^n x_iy_i$  and $A\!:\!B:= {\rm Tr}(A B^\top)$. For a function $f : E \to F$, where $E$ and $F$ are arbitrary vector spaces, we define: 
 \begin{equation}\label{notation} 
 \big\{ f \big\}_x^y := f(y) - f(x) \ \mbox{for all} \ x,y \in E. 
 \end{equation}
 Throughout all the paper, partitions of an interval $[0,T]$ are generically denoted by $\pi = \{t_0, \dots, t_N\}$, the mesh of $\pi$ by $| \pi | := \sup_{0 \le n \le N-1} | t_{n+1} - t_n |$, and we denote the closest left point of the grid to any point $s\in[0,T]$ by $t_{n-1}(s):=\sup\{t \in\pi:t\le  s\}$.  


\section{Itô's formula}\label{sec:Ito}

\begin{definition}
A map $u : \cP_2(\dbR^d) \longrightarrow \dbR$ has a linear functional derivative $\d_m u : \cP_2(\dbR^d) \times \dbR^d \to \dbR$ if $\d_m u$ is quadratically growing in $x$, locally uniformly in $\mu$, and
$$ u(\m_0) - u(\m_1) 
= \int_0^1 \int_\dbR \d_m u(\m_\l, x) (\m_0\!-\!\m_1)(dx)d\l \q \mbox{for all} \ \m_0, \m_1 \in \cP_2(\dbR^d), 
$$
with $\m_\l := (1\!-\!\l)\m_0 + \l \m_1$.
\end{definition}
Using the notation \eqref{notation}, we may rewrite this through arbitrary r.v. $\xi_0,\xi_1$ with distribution $\mu_0$ and $\mu_1$, respectively:
\be\label{deltamu}
\{u\}_{\mu_1}^{\mu_0}
=
\int_0^1 \dbE\big[\big\{\d_m u(\m_\l, \cdot)\big\}_{\xi_1}^{\xi_0}\big]d\l.
\ee
Similarly, we define the map $(m,x,\hat x)\in\cP_2(\dbR^d) \times \dbR^d \times\dbR^d \longmapsto \delta^2_{mm}u(m,x,\hat x)$ as the linear functional derivative of $m\longmapsto\d_m u(m,x)$ with respect to $m$. We denote by $C_b^2(\cP_2(\dbR^d))$ the set of functions $u : \cP_2(\dbR^d) \to \dbR$ such that $\d_m u$, $\pa_x u$, $\pa_{xx}^2 \d_m u$ and $\pa_{x\hat x}^2 \d_{mm}^2 u$ exist and are continuous and bounded. 

Given a stochastic process $X : [0,T] \times \O \longrightarrow \dbR^d$ and a random probability measure map $m: [0,T]\times \O\longrightarrow\cP_2(\dbR^d)$, we denote
$$
\Delta X_s := X_s - X_{s-}=\{X\}_{s-}^s,~~
\Delta m_s := \{m\}_{s-}^s,
~~\mbox{and}~~
J_t(\bm{m}):=\{s\le t: \Delta m_s\neq 0\}.
$$
Moreover, given a map $u : \cP_2(\dbR^d) \longrightarrow \dbR$, we also denote:
$$
\Delta u(m_s) :=\{u\}^{m_s}_{m_{s-}}, 
~~\mbox{and}~~
\D^{\! X} \d_m u(m_s, X_s) := \big\{\d_m u(m_s, \cdot)\big\}^{X_s}_{X_{s-}}. \\ 
$$

\begin{theorem}\label{thm:Ito}
Let $X \!:=\! V \!+\! M$ be a càdlàg semimartingale with finite variation process $V$ and local martingale $M$ satisfying:
\bea\label{integrability} 
\dbE\Big[ | V_T^c |^2 + {\rm Tr}[M]_T^c + \Big(\sum_{0 < s \le T} |\D X_s| \Big)^2\Big] < \infty.
\eea
Given a sub-sigma algebra $\cF^0 \subset \cF$, denote by $m_t$ the conditional law of $X_t$ given $\cF^0$, $\dbE^0[ \cdot ] := \dbE^0[ \cdot | \cF_T^0]$ and $\hat \dbE^0[ \cdot ] := \dbE^0[ \cdot | \cF_T^0, X]$. 
Let $u \in C_b^2(\cP_2(\dbR^d))$ be such that $\pa_x \pa_{\hat x} \d_{mm}^2 u(\cdot, x, \hat x) \in C^0(\cP_1(\dbR^d))$ for all $x, \hat x \in \dbR^d$. Then, we have:
\begin{align*}
\{u\}^{m_t}_{m_0} 
= \dbL_t^X u + \frac12 \dbE^0 \hat \dbE^0\Big[ \int_0^t \pa^2_{x\hat x} \d_{mm}^2 u(m_s, X_s,\hat X)\!:\! d[ X, \hat X ]^c_s 
\Big], \ \mbox{a.s.  for all}~ t \in [0,T],
\end{align*}
where $\hat X$ denotes an independent copy of $X$ conditionally on $\cF^0_T$, and:
\begin{eqnarray*}
\dbL^X_t u &:=&\sum_{s\le t} \D u(m_s)
                      +\dbE^0\big[L_t^X \d_m u \big],
\\
L_t^X \f &:=&\sum_{s \in J_t(\bm{m})^c} \D^{\!X}\! \f (m_s, \!X_s)
                         + \int_0^t \pa_x \f(m_s, \!X_s) \!\cdot\! dX_s^c 
                                     + \frac12\pa_{xx}^2\f(m_s, \!X_s) 
                                                   \!:\! d[X]_s^c. 
\end{eqnarray*}
\end{theorem}
The proof of this result uses the following additional notations. Given a partition $\pi=\{t_0,\ldots, t_N\}$, we denote for $r\in[0,1]$:
$$
\Delta^{\!\pi}X_{t_n}\!:=\!X_{t_n}\!\!-X_{t_{n-1}},
~X_{t_{n-1},t_n}^{r}\!:=\!X_{t_{n-1}}\!+\!r\Delta^{\!\pi}X_{t_n},
~X_s^r\!:=\!X_{s-}\!+\!r\Delta X_s.
$$
We denote similarly $\Delta^{\!\pi}\hat X_{t_n},
~\hat X_{t_{n-1},t_n}^{r},~\hat X_s^r$, and:
$$
\Delta^{\!\pi}m_{t_n}\!:=\!m_{t_n}\!\!-m_{t_{n-1}},
~m_{t_{n-1},t_n}^{r}\!:=\!m_{t_{n-1}}\!+\! r\Delta^{\!\pi}m_{t_n},
~\mbox{and}~
m_s^r\!:=\!m_{s-}\!+\!r\Delta m_s
$$
The following Lemma (proved in Appendix \ref{sec:appendix}) play a crucial role in the proof of Theorem \ref{thm:Ito}.

\begin{lemma}\label{lem:cv-bracket}
Let $X$ and $\hat X$ be as in Theorem \ref{thm:Ito}, and consider a dense partition $\pi = \{t_0, \dots, t_N\}$ of $[0,T]$.
\\
{\rm (i)} Let $H$ be a bounded progressively measurable left-limited process and $\{ H^\pi \}_\pi$ a uniformly bounded family of simple adapted processes with $H_{t_{n-1}(s)}^\pi \longrightarrow H_{s-}$ as $| \pi | \to 0$, a.s. Then:
$$ \sum_{n=1}^N H_{t_{n-1}}^\pi\!:\! \D^{\!\pi}\!X_{t_n}(\D^{\!\pi}\!X_{t_n})^\intercal   
\longrightarrow \int_0^T H_{s-}\!:\!d[X]_s, \q \mbox{as}~| \pi | \to 0~\mbox{in $\dbL^1$}.  $$ 
{\rm (ii)} Let $f: \cP_2(\dbR^d) \times \dbR^d \times \dbR^d\longrightarrow\dbR$ be a bounded continuous map such that $f \in C^0(\cP_1(\dbR^d) \times \dbR^d \times \dbR^d)$, and denote $\theta_n^r:=(m_{t_{n-1}, t_n}^{r_1}, X_{t_{n-1}, t_n}^{r_2}, \hat X_{t_{n-1}, t_n}^{r_3})$, $\theta_s^r:=(m_s^{r_1}, X_s^{r_2}, \hat X_s^{r_3})$ for $r=(r_1,r_2,r_3)\in[0,1]^3$. Then, there exists a subsequence of partitions, still named $\pi$, such that
$$ 
\sum_{n=1}^{N} \big\{f\}^{\theta_n^r}_{\theta_n^0}\,\D^{\!\pi} X_{t_n} \!\cdot\! \D^{\!\pi} \hat X_{t_n} 
\longrightarrow
\sum_{0 < s \le T} \big\{f\}^{\theta_s^r}_{\theta_s^0}\, \D X_s \!\cdot\!\D\hat X_s,
\q \mbox{as}~| \pi | \to 0~\mbox{in}~\dbL^1. 
$$
\end{lemma} 

\noindent\textbf{Proof of Theorem \ref{thm:Ito}.}
For simplicity, we prove the formula in the case $d=1$ as the only difficulty in higher dimension is related to notations. Let $\pi := \{t_0, \dots, t_N\}$ be a partition of $[0,T]$.  Fix $n \in [N-1]$. By definition of the functional linear derivative in \eqref{deltamu}, we have:
\begin{equation}\label{decomposition}
\{u\}^{m_{t_N}}_{m_{t_0}}
=
\sum_{n=1}^N \{u\}^{m_{t_n}}_{m_{t_{n-1}}} 
=
\sum_{n=1}^N \int_0^1 \dbE^0\Big[ \big\{ \d_m u(m_{t_{n-1}, t_n}^{\l_1}, \cdot) \big\}_{X_{t_{n-1}}}^{X_{t_n}} \Big]d\l_1 
= 
\sum_{n=1}^N 
S_n+T_n^0 - T_n^1,
\end{equation}
with
$$
S_n\!:=\!\!\int_0^1\!\dbE^0\Big[ \big\{ \d_m u(m_{t_{n-1}}, \cdot) \big\}_{X_{t_{n-1}}}^{X_{t_n}}\Big]d\l_1,
~\mbox{and}~
T_n^j \!=\!\!\int_0^1 \!\!\dbE^0\Big[\big\{ \d_m u(\cdot, X_{t_{n-j}}) \big\}_{m_{t_{n-1}}}^{m_{t_{n-1}, t_n}^{\l_1}}\Big]d\l_1,~j\!=\!0,1.
$$
{\bf 1.} By the standard Itô formula inside the conditional expectation, we obtain for the first term of \eqref{decomposition}:
\begin{align*}
S_n \!=\! 
\dbE^0 \Big[\int_{t_{n-1}}^{t_n} \!\!\pa_x \d_m u(m_{t_{n-1}}, X_s)dX_s^c 
\!+\! \frac 1 2 \pa_{xx}^2 \d_m u(m_{t_{n-1}}, X_s)d[X]_s^c 
\!+ \!\!\!\!\sum_{(t_{n-1},t_n]} \!\!\!\D^{\!X} \d_m u(m_{t_{n-1}}, X_s) \Big].
\end{align*}
Then, in view of our growth conditions on $u$ and its derivatives, we deduce from the dominated convergence theorem that 
\begin{align}\label{limSn}
\sum_{n=1}^{N} S_n 
\underset{| \pi | \to 0}{\longrightarrow}  
\dbE^0 \Big[\int_0^T \!\!\pa_x \d_m u(m_s, X_s)dX_s^c 
                                \!+\! \frac 1 2 \pa_{xx}^2 \d_m u(m_s, X_s)d[X]_s^c 
                                \!+\! \sum_{(0,T]} \D^X \d_m u(m_{s-}, X_s) \Big].
\end{align}
{\bf 2.} We next focus on the remaining two terms in \eqref{decomposition}. Using again the definition of the functional linear derivative, together with standard finite dimensional differential calculus, we obtain for $j=0,1$ that:
\begin{align*}
T_n^j 
&= \int_{[0,1]^2}  \dbE^0 \hat \dbE^0\Big[\big\{ \d^2_{mm} u(m_{t_{n-1}, t_n}^{\l_1 \l_2}, X_{t_{n-j}}, \cdot) \big\}_{\hat X_{t_{n-1}}}^{\hat X_{t_n}} \Big]\l_1d\l_1 d\l_2 \\
&= \int_{[0,1]^3} \dbE^0 \hat \dbE^0\big[ \pa_{\hat x}\d^2_{mm} u(m_{t_{n-1}, t_n}^{\l_1 \l_2}, X_{t_{n-j}}, \hat X_{t_{n-1}, t_n}^{\l_3}))\D^{\!\pi} \hat X_{t_n}\big]\l_1d\l_1 d\l_2 d\l_3.
\end{align*} 
Then, it follows again from standard finite dimensional differential calculus that:
\begin{eqnarray} 
T_n^0 - T_n^1 
&=& 
\int_{[0,1]^4}  \dbE^0 \hat \dbE^0\big[ f(m_{t_{n-1}, t_n}^{\l_1 \l_2}, X_{t_{n-1}, t_n}^{\l_4}, \hat X_{t_{n-1}, t_n}^{\l_3})) \D^{\!\pi} X_{t_n} \D^{\!\pi} \hat X_{t_n}\big]  \l_1 d\bm{\l},
~
f:=\pa_{x}\pa_{\hat x}\d^2_{mm} u,
\nonumber\\
&=& 
\int_{[0,1]^4}  \dbE^0 \hat \dbE^0\big[ f(\th_n^0) \D^{\!\pi} X_{t_n} \D^{\!\pi} \hat X_{t_n} + R_n^\pi \big]  \l_1 d\bm{\l}
\label{Tbarn}
\end{eqnarray}
with $\bm{\l} := (\l_1, \l_2, \l_3, \l_4)$ and using the notations of Lemma \ref{lem:cv-bracket} (ii):
$$ 
R_n^\pi
:= 
\big\{f\big\}^{\theta^r_n}_{\theta^0_n}\,\,\D^{\!\pi} X_{t_n} \!\cdot\! \D^{\!\pi} \hat X_{t_n}
,~~r:=(\l_1\l_2,\l_4,\l_3). 
$$
{\bf 3.} As $f = \pa_x \pa_{\hat x} \d_{mm}^2 u$ satisfies the conditions of Lemma \ref{lem:cv-bracket} (ii), we deduce that, after possibly passing to a subsequence, the following convergence holds in $\dbL^1$, and therefore under conditional expectation:
\begin{align*}
\sum_{n=1}^{N} f(\theta^0_{n}) \D^{\!\pi}X_{t_n} \D^{\!\pi}\hat X_{t_n} + R_n^\pi 
\underset{| \pi | \to 0}{\longrightarrow} 
&\int_0^T f(\theta_s^0) d[X, \hat X]_s 
+ \sum_{0 < s \le T} \big\{f\}^{\theta_s^r}_{\theta_s^0}\, \D X_s \D\hat X_s
\\
&= 
\int_0^T f(\theta^{0}_s) d[X, \hat X]_s^c 
+ \sum_{0 < s \le T} f(\theta_s^r)\D X_s \D \hat X_s,
\end{align*}
As $f$ is bounded, it follows from the dominated convergence theorem that:
\begin{align}\label{limbarTn}
\sum_{n=1}^N T_n^0 - T_n^1
\underset{| \pi | \to 0}{\longrightarrow} \frac12 \dbE^0 \hat \dbE^0\Big[ \int_0^T f(\th_s^0)d[X,\hat X]_s^c \Big] +
\int_{[0,1]^4} \dbE^0 \hat \dbE^0\big[\sum_{0 < s \le T} f(\theta_s^r)\D X_s \D \hat X_s \Big]\l_1 d\bm{\l}.
\end{align}
Recall that $f=\pa_{x}\pa_{\hat x}\d^2_{mm} u$, and observe that 
\begin{align}
\int_{[0,1]^2}  f(\theta_s^r)\D X_s \D \hat X_s d\l_3 d\l_4 
&= \big\{ \d_{mm}^2 u(m_s^{\l_1 \l_2}, X_s, \cdot) \big\}_{\hat X_{s-}}^{\hat X_s} 
    - \big\{  \d_{mm}^2 u(m_s^{\l_1 \l_2}, X_{s-}, \cdot )  \big\}_{\hat X_{s-}}^{\hat X_s}
\label{jumpsXhatX}\\
&\hspace{-10mm}
   = \1_{\{s\in J(m)\}}\Big(\big\{ \d_{mm}^2 u(m_s^{\l_1 \l_2}, X_s, \cdot) \big\}_{\hat X_{s-}}^{\hat X_s} 
    - \big\{  \d_{mm}^2 u(m_s^{\l_1 \l_2}, X_{s-}, \cdot )  \big\}_{\hat X_{s-}}^{\hat X_s}\Big),
\nonumber
\end{align}
as $m_s^\l = m_s$ whenever $s \in J(\bm{m})^c$, and each of these terms are nonzero if and only if $X$ and $\hat X$ jump at the same time, which happens with probability zero on $J(\bm{m})^c$ since $X$ and $\hat X$ are independent conditionally on $\cF^0$.  
Note also that, as the jumps of the flow $s \mapsto m_s$ are $\cF_T^0$-measurable, 
\begin{align*}
&\int_{[0,1]^2} \dbE^0 \hat \dbE^0\Big[ \sum_{s \in J(\bm{m})} \big\{ \d_{mm}^2 u(m_s^{\l_1 \l_2}, X_s, \cdot) \big\}_{\hat X_{s-}}^{\hat X_s} 
- \big\{  \d_{mm}^2 u(m_s^{\l_1 \l_2}, X_{s-}, \cdot )  \big\}_{\hat X_{s-}}^{\hat X_s}\Big]\l_1 d\l_1 d\l_2 \\ 
&= \sum_{s \in J(\bm{m})} \int_{[0,1]^2} \dbE^0 \hat \dbE^0\Big[  \big\{ \d_{mm}^2 u(m_s^{\l_1 \l_2}, X_s, \cdot) \big\}_{\hat X_{s-}}^{\hat X_s} 
- \big\{  \d_{mm}^2 u(m_s^{\l_1 \l_2}, X_{s-}, \cdot )  \big\}_{\hat X_{s-}}^{\hat X_s}\Big]\l_1 d\l_1 d\l_2 \\
&=  \sum_{s \in J(\bm{m})} \dbE^0\Big[ \big\{ \d_m u(\cdot, X_s)\big\}^{m_s^{\l_1}}_{m_{s-}} - \big\{ \d_m u(\cdot, X_{s-})\big\}^{m_s^{\l_1}}_{m_{s-}}\Big]d\l_1 \\
&= \sum_{s \in J(\bm{m})} \D u(m_s)  - \dbE^0 \big[ \D^X \d_m u(m_{s-}, X_s) \big].
\end{align*}
Therefore, we have:
\bea\label{Tn12} 
\sum_{n=1}^N  T_n^0 - T_n^1
&\underset{| \pi | \to 0}{\longrightarrow}& 
\sum_{s \in J(\bm{m})} \D u(m_s) - \dbE^0 \big[\D^X \delta_m u(m_{s-}, X_s)\big], \nonumber
\eea
Together with \eqref{limSn}, \eqref{Tbarn} and \eqref{limbarTn}, this provides the required result by taking limits in \eqref{decomposition} along appropriate subsequences. 
\qed

\begin{example}\label{ex:poisson} {\rm
Let $X$ be an $\dbR$-valued jump-diffusion satisfying the following dynamics:
$$ dX_t = b_t dt + \si_t dW_t + \g_t dN_t + \si_t^0 dW_t^0 + \g_t^0 dN_t^0, $$
where $W, W^0$ are two independent Brownian motions in $\dbR$, and $N, N^0$ are two independent pure jump processes in $\dbR$ with bounded intensities $\l$ and $\l^0$ from $\dbR_+\times \O$ to $\dbR$, with jump size following distributions $\n$ and $\n^0$, assumed to have finite second order moments. For simplicity, we assume that $b, \si, \si^0, \g$ and $\g^0$ are bounded adapted processes taking values in $\dbR$. Let $\dbF^0$ be the filtration generated by $(W^0, N^0)$, and $u$ a smooth map from $\cP_2(\dbR)$ to $\dbR$. We denote by $m_t$ the law of $X_t$ conditionally on $\dbF^0$.  

We first verify that $X$ satisfies the integrability conditions \eqref{integrability}. It is clear that:
$$ \dbE\Big[ | V_T^c |^2 + [X]_T^c \Big] < \infty.  $$
Also observe that:
$$ \Big( \sum_{0 < s \le T} \D X_s \Big)^2 \le 2 \Big( \sum_{0 < s \le T} \g_{s-} \D N_s \Big)^2 +2 \Big( \sum_{0 < s \le T} \g_{s-}^0 \D N_s^0 \Big)^2. $$
Let $\Pi$ be the random Poisson measure associated with $N$, that is such that $N_t = \int_0^t \int_{\dbR \times \dbR_+} y \1_{\th \le \l_s} \Pi(ds, d\th, dy)$. We assume that its compensator writes $ds d\th \n(dy)$. Then, introducing $\bar \Pi(ds,d\th,dy) :=  \Pi(ds,d\th,dy) - ds d\th \n(dy)$ the compensated measure of $\Pi$, we have the following estimates:
\begin{eqnarray*}
\dbE\Big[ \Big( \sum_{0 < s \le T} \g_{s-} \D N_s \Big)^2  \Big] 
&\hspace{-2mm}=&\hspace{-2mm}
 \dbE\Big[ \Big( \int_0^T \int_{\dbR_+ \times \dbR} \g_{s-} y \1_{\th \le \lambda_s} \Pi(ds, d\th, dy)  \Big)^2  \Big] 
\\
&\hspace{-2mm}\le&\hspace{-2mm} 
\dbE\Big[ 2 \int_0^T\!\!\!\int_{\dbR_+ \times \dbR} \g_{s-} y \1_{\th \le \lambda_s} \bar \Pi(ds, d\th, dy)  \Big)^2 \!\!+\!  2T \int_0^T\!\!\!\int_\dbR \g_{s-}^2 y^2 \l_s^2 ds \n(dy)  \Big] 
\\
&\hspace{-2mm}=&\hspace{-2mm}
  \dbE\Big[ 2 \int_0^T\!\!\!\int_\dbR \g_{s-}^2 y^2 \l_s ds \n(dy)  +  2T \int_0^T\!\!\!\int_\dbR \g_{s-}^2 y^2 \l_s^2 ds \n(dy)  \Big] < \infty,
\end{eqnarray*}
by boundedness of $\g$ and $\l$ and integrability of $\n$. We have the same estimate for $N^0$, and therefore \eqref{integrability} is satisfied.

Then, according to Theorem \ref{thm:Ito}, we have:
\begin{align*}
u(m_t) =& u(m_0) + \dbL_t^X u 
                 + \frac12 \dbE^0 \hat \dbE^0\Big[ \int_0^t \pa^2_{x\hat x} 
                                                                      \d_{mm}^2 u(m_s, X_s,\hat X)\!:\! 
                                                                      \sigma^0_s(\hat\sigma^0_s)^\intercal ds
                                                              \Big], \q \dbP-\mbox{a.s.,}
\end{align*}
with
\begin{eqnarray*}
&&\dbL^X_t u 
:= \sum_{s \le t} \D u(m_s)
                      +\dbE^0\big[L_t^X \d_m u \big],
\\
&&L_t^X \!\f 
:= \!\!\sum_{s \in J_t(\bm{m})^c} \D^{\!X} \!\f (m_s, \!X_s)
                         +\int_0^t \cL_s^c\varphi(m_s,X_s)ds
                         +\pa_x \f(m_s, \!X_s) \!\cdot\! \sigma^0_sdW^0_s 
\\
&&\cL_s^c\varphi(m_s,X_s)
:= b_s\!\cdot\!\partial_x\varphi(m_s,X_s)
                                                + \frac12(\sigma_s\sigma_s^{\intercal}
                                                               +\sigma^0_s{\sigma^0_s}^{\intercal})
                                                             \!:\!
                                                             \pa_{xx}^2\f(m_s, \!X_s).
\end{eqnarray*}                                                   
}
\end{example}

\section{Itô-Wentzell's formula}\label{sec:Ito-Wentzell}

In this paragraph, we provides an extension of Theorem \ref{thm:Ito} to the case where the function $u$ is a random field.

 More precisely, we consider $U : [0,T] \times \cP_2(\dbR^d) \times \O \to \dbR$ such that:
$$ U_t(m) = U_0(m) + \int_0^t \phi_{s}(m)  dA_s + \int_0^t \psi_{s}(m)  dN_s, $$
with $\phi, \psi : [0,T] \times \cP_2(\dbR^d) \to \dbR$ two predictable maps, $A : [0,T] \times \O \to \dbR$ a finite variation process and $N : [0,T] \times \O \to \dbR$ a local martingale.

\begin{assumption}\label{assum:ItoWentzell}
The maps $f \in \{U_0, \phi_t, \psi_t, t \in [0,T]\}$, defined on $\cP_2(\dbR^d)$, satisfy:\\
{\rm (i)} $\d_m f$, $\pa_{xx}^2 \d_m f$, $\d_{mm}^2 f$ and $\pa_{x\hat x}^2 \d_{mm}^2 f$ exist and are continuous, \\
{\rm (ii)} $f$, $\pa_{xx}^2 \d_m f$ and $\pa_{x\hat x}^2 \d_{mm}^2 f$ are bounded. \\
Moreover, the random field $U$ is such that: \\
{\rm (iii)} $| A |_{\rm{TV}}$ and $[N]_T$ are bounded, almost surely. \\
\end{assumption}

Similarly to \citeauthor*{fadle2024ito} \cite{fadle2024ito}, Assumption \ref{assum:ItoWentzell} implies that $U$ is differentiable with respect to $m$, and that the following hold: 
\begin{eqnarray*}
\pa_{x}^i \d_m U_t(m,x) 
&\hspace{-2mm}=& \hspace{-2mm} 
\int_0^t \pa_{x}^i\d_m \phi_{s}(m,x)  dA_s 
    + \int_0^t \pa_{x}^i\d_m \psi_{s}(m,x) dN_s, ~ i=0,1,2, \\
\pa_{x}^i \pa_{{\hat x}}^i \d_{mm}^2 U_t(m,x, \hat x) 
&\hspace{-2mm}=&\hspace{-2mm}
\int_0^t \pa_{x}^i\pa_{{\hat x}}^i \d_{mm}^2 \phi_{s}(m,x, \hat x)  dA_s 
    + \int_0^t \pa_{x}^i \pa_{{\hat x}}^i\d_{mm}^2 \psi_{s}(m,x, \hat x)  dN_s, ~ i=0,1,
\end{eqnarray*}
and all these functions are continuous. Here $\pa^0_{x} \d_m:=\d_m$, $\pa^1_{x} \d_m:=\pa_{x} \d_m$, $\pa^2_{x} \d_m:=\partial^2_{xx}\delta_m$, and similar notations for $\pa_{x} \pa^i_{{\hat x}} \d_{mm}^2 $, $i=0,1$. 

We now introduce the following definition which is implicitly adopted throughout the existing literature on this topic.

\begin{definition}
 Let $F$ be a random field defined on some measurable space $(E, \cA)$, and let $\xi \in L^2(\O, \cA)$. We define $F_t(\xi) := F_t(x) \big|_{x = \xi}$, i.e. $F_t(\xi)(\omega):= F_t(\xi(\omega))(\omega)$ for $\dbP-$a.e. $\omega\in\Omega$. 
 \end{definition}

\begin{remark}\label{rem:expectation-random-field}{\rm Assume for simplicity that $F$ is a stochastic integral $F_t(x) = \int_0^t g_s(x)dW_s$, for all $x \in E$, for some Brownian motion $W$. Then, given that the process $\{g_s(\xi),s\le t\}$ is anticipative, $F_t(\xi) = \int_0^t g_s(\xi)\, dW_s$ is \textbf{not} defined in the It\^o stochastic integration sense, and it is not meant either in the Skorohod sense. Moreover, we notice that the last definition implies that \(F_t(\xi)\) is \(\mathcal{F}_t\)-measurable whenever $\xi$ is $\cF_t$-measurable.
   
}
\end{remark}

The following Lemma, proved in the Appendix \ref{sec:appendix}, will be crucial in the proof of Theorem \ref{thm:Ito-Wentzell}:
\begin{lemma}\label{lem:cv2}
{\rm (i)} Let $F : [0,T] \times \cP_2(\dbR^d) \times \dbR^d \times \O \longrightarrow \dbR$ be a random field defined through continuous bounded  predictable processes $f, g$, and some finite variation process $A$ and a local martingale $N$ satisfying Assumption \ref{assum:ItoWentzell} and independent of $X$ conditionally on $\cF^0$ by:
$$
  F_t(\cdot) = \int_0^t f_{s}(\cdot)  dA_s + \int_0^t g_{s}(\cdot)  dN_s.
$$
Then, for any dense partition $\pi$ and $\l \in [0,1]$, we have $\dbP-$a.s.
\begin{align*}
\dbE^0\Big[\sum_{n=1}^{N} \big\{F_.(m_{t_{n-1}},\cdot)\big\}^{t_n,X_{t_{n-1},t_n}^{\l}}_{t_{n-1},X_{t_{n-1}}}
| \D^{\!\pi}\!X_{t_n}|^2\Big] 
\underset{| \pi | \searrow 0} {\longrightarrow} 
\dbE^0\Big[\sum_{0 < s \le T} \big\{F_.(m_{s-}, \cdot)\big\}^{s,X_{s}^{\l}}_{s-,X_{s-}}
| \D X_{s}|^2\Big].
\end{align*}
{\rm(ii)}  Let $\hat F : [0,T] \times \cP_2(\dbR^d) \times \dbR^d \times \dbR^d \times \O \longrightarrow \dbR$ be a random field defined through continuous bounded  predictable processes $\hat f,\hat g$, and some finite variation process $A$ and a local martingale $N$ satisfying Assumption \ref{assum:ItoWentzell} and independent of $(X, \hat X)$ conditionally on $\cF^0$:
$$
  \hat F_t(\cdot) = \int_0^t \hat f_{s}(\cdot)  dA_s + \int_0^t \hat g_{s}(\cdot)  dN_s.
$$
We furthermore assume that $\hat f$ and $\hat g$ are continuous with respect to the 1-Wasserstein distance $\cW_1$. Then, for any dense partition $\pi$ and $r := (r_1, r_2, r_3) \in [0,1]^3$, we have $\dbP-$a.s.
\begin{align*}
\dbE^0 \hat \dbE^0 
\Big[\sum_{n=1}^{N} \big\{\hat F_{.}(\cdot)\big\}^{t_n,\th_n^r}
                                                                                       _{t_{n-1},\th_n^{\bm{0}}}
                                  \D^{\!\pi}\!X_{t_n}\D^{\!\pi}\!\hat X_{t_n}
\Big] 
\underset{| \pi | \searrow 0} {\longrightarrow} 
\dbE^0 \hat \dbE^0
\Big[\sum_{0 < s \le T} \big\{\hat F_{.}(\cdot)\}^{s,\th_s^r}
                                                                                                _{s-,\th_s^{\bm{0}}}
                                    \D X_s\D\hat X_s
\Big],
\end{align*}
where $\th_n^r := (m_{t_{n-1}, t_n}^{r_1}, X_{t_{n-1},t_n}^{r_2},\hat X_{t_{n-1},t_n} ^{r_3})$ and $\th_s^r := (m_{s}^{r_1}, X_{s}^{r_2},\hat X_{s}^{r_3})$. 
\end{lemma} 
 
\begin{theorem}\label{thm:Ito-Wentzell}
Let $X$ be a càdlàg semimartingale satisfying the integrability conditions \eqref{integrability}, and let $m_t$ denote the conditional law of $X$ given $\cF^0$. Further assume that, for all compact $K \subset \cP_2(\dbR^d)$,
\bea\label{eq:finiteness-derivatives}
\sup_{ (t,m) \in [0,T] \times K} \;\sup_{x, \hat x \in \dbR^d} \Big\{ | \pa_x \d_m U_t | +  | \pa_{xx}^2 \d_m U_t | + | \pa_{x \hat x}^2 \d_{mm}^2 U_t | \Big\}(m,x, \hat x) < \infty, \ \dbP\mbox{-a.s}.
\eea
Under Assumption \ref{assum:ItoWentzell}, we have:
\begin{eqnarray*}
\{U\}^{t,m_t}_{0,m_0} 
&\hspace{-2mm}=&\hspace{-2mm}
\dbL_t^X U + \frac12 \dbE^0 \hat \dbE^0\Big[ \int_0^t \pa^2_{x\hat x} \d_{mm}^2U_s
                                                                    (m_s,  X_s)\!:\!d[  X, \hat X ]^c_s 
                                               \Big]  
\\ 
&\hspace{-2mm}&\hspace{-2mm}
+ \int_0^t \phi_{s}(m_{s-})\!\cdot\! dA_s + \psi_{s}(m_{s-})\!\cdot\!dN_s
\\
&\hspace{-2mm}&\hspace{-2mm}
+\dbE^0\Big[\sum_{0 < s \le T} \!\!\pa_x \d_m \phi_{s}(m_{s-}, X_{s-}) 
                                                      \!\cdot\! \D  X_s \D A_s
                         \!+\!\! \int_0^t\!\! \pa_x \d_m \psi_{s}(m_{s-},  X_{s-})\!:\!d[ X, N]_s
                  \Big], 
~\mbox{a.s.,}
\end{eqnarray*}
where $\dbE^0[ \cdot ] := \dbE[ \cdot | \cF_T^0, A, N]$, $\hat \dbE^0[ \cdot ] := \dbE[ \cdot | \cF_T^0, A, N, X]$ and $\hat X$ denotes a copy of $X$, independent from $X$ conditionally on $\cF^0_T$. 

\proof 
 For an arbitrary dense partition $\pi := \{t_0, \dots, t_N\}$, we have:
\begin{eqnarray}\label{eq:decomposition}
 \!\!\!\! \!\!\!\!
 \{U\}^{t_N,m_{t_N}}_{t_0,m_{t_0}} 
 = \sum_{n=1}^{N}\big\{ U_\cdot (m_{t_{n-1}}) \big\}_{t_{n-1}}^{t_n} 
       \!+\!  \big\{ U_{t_n}\}_{m_{t_{n-1}}}^{m_{t_n}} 
 &\!\!\!\!\!=&\!\!\!\!\! \int_0^T dU_s(m_{t_{n-1}(s)}) 
    + \sum_{n=1}^{N}\big\{ U_{t_n}\}_{m_{t_{n-1}}}^{m_{t_n}}
\nonumber\\
&\hspace{-68mm}=&\hspace{-35mm}
\int_0^T \!\!\!\phi_{s}(m_{t_{n-1}(s)})dA_s \!+\! \psi_{s}(m_{t_{n-1}(s)})dN_s
     +\! \sum_{n=1}^{N}\!\big\{ U_{t_n}\}_{m_{t_{n-1}}}^{m_{t_n}}. 
\end{eqnarray}
Since $\cW_2(m_{t_{n-1}(s)}, m_{s-}) \to 0$ as $| \pi | \to 0$, a.s., we deduce from the boundedness of $\phi$ and the dominated convergence theorem for Stieltjes integrals that $\int_0^T \phi_{s}(m_{t_{n-1}(s)})dA_s \to \int_0^T \phi_{s}(m_{s-})dA_s$, a.s. As for the stochastic integral, we see that:
$$ \dbE\Big[\Big(\int_0^T \psi_{s}(m_{t_{n-1}(s)})dN_s - \int_0^T \psi_{s}(m_{s-})dN_s\Big)^2\Big] \le \dbE\Big[ \int_0^T \big| \psi_{s}(m_{t_{n-1}(s)}) - \psi_{s}(m_{s-}) \big|^2 d[ N ]_s \Big],$$
 and we then deduce from the dominated convergence theorem for Stieltjes integrals that $\int_0^T \psi_{s}(m_{t_{n-1}(s)})dN_s \to  \int_0^T \psi_{s}(m_{s-})dN_s$ in $\dbL^2$, and therefore a.s. along some subsequence (still denoted the same). 
 
Denoting $G:=\pa_x \d_m U$, we rewrite the last term on the right hand side of \eqref{eq:decomposition} as:
\bea\label{terms}
\big\{ U_{t_n}\}_{m_{t_{n-1}}}^{m_{t_n}} 
&=&
\int_{[0,1]^2} \dbE^0\Big[ G_{t_n}(m_{t_{n-1}, t_n}^{\l_1}, X_{t_{n-1},t_n}^{\l_2})\D^{\!\pi} X_{t_n}\Big]d\l_1 d\l_2 \nonumber \\
&=& 
\dbE^0\Big[ G_{t_n}(m_{t_{n-1}}, X_{t_{n-1}})\D^{\!\pi} \!X_{t_n}\Big]
\nonumber\\
&&+  \int_{[0,1]^2}\!\! \dbE^0\Big[ \{ G_{t_n}(\cdot, X_{t_{n-1},t_n}^{\l_2}) \}_{m_{t_{n-1}}}^{m_{t_{n-1}, t_n}^{\l_1}} \D^{\!\pi} \!X_{t_n}\Big]d\l_1 d\l_2 
+ \dbE^0[R_n],
\eea
where, denoting $F:=\partial_xG=\pa_{xx}^2 \d_m U$,
\bea
R_n
&=&
\int_{[0,1]^2} \dbE^0\Big[ \big\{ G_{t_n}(m_{t_{n-1}}, \cdot) 
                                            \big\}_{X_{t_{n-1}}}^{X_{t_{n-1},t_n}^{\l_2}} \D^{\!\pi}\!X_{t_n}
                                    \Big]d\l_1 d\l_2 
\nonumber \\
&=&
\int_{[0,1]^3} \dbE^0\Big[ F_{t_n}(m_{t_{n-1}}, X_{t_{n-1}, t_n}^{\l_2 \l_3})|
                                        \D^{\!\pi}\!X_{t_n}|^2\Big]\l_3 d\l_1 d\l_2 d\l_3
\nonumber \\
&=&
\frac{1}{2} \dbE^0\Big[ F_{t_{n-1}}(m_{t_{n-1}}, X_{t_{n-1}})| \D^{\!\pi} \!X_{t_n}|^2\Big] 
\nonumber \\
&&
+\! \int_{[0,1]^3} \!\!\dbE^0\Big[ \big\{F_{.}(m_{t_{n-1}},.)\big\}^{t_n,X_{t_{n-1}, t_n}^{\l_2 \l_3}}                                                                                                                             _{t_{n-1},X_{t_{n-1}}}
                                           | \D^{\!\pi} \!X_{t_n}|^2
                                   \Big]\l_3 d\l_2 d\l_3.
\nonumber
\eea
 \begin{enumerate}
 \item Decompose the first term in \eqref{terms} as:
 \begin{eqnarray*}
 \dbE^0\Big[ G_{t_n}(m_{t_{n-1}}, X_{t_{n-1}})\D^{\!\pi}\!X_{t_n}\Big] 
 &=& 
 \dbE^0\Big[ G_{t_{n-1}}(m_{t_{n-1}}, X_{t_{n-1}})\D^{\!\pi}\!X_{t_n}\Big]
 \\
 &&+ \dbE^0\Big[\big\{ G_{\cdot}(m_{t_{n-1}}, X_{t_{n-1}})
                          \big\}_{t_{n-1}}^{t_n}\D^{\!\pi}\!X_{t_n}\Big].
\end{eqnarray*}
Using \eqref{eq:finiteness-derivatives} with $K$ the closure of the graph of $s \mapsto m_s$ which is compact in $\cP_2(\dbR^d)$, a.s, we have $ \dbE^0\Big[ \int_0^T | G_s(m_s, X_s) |^2 ds \Big] < \infty$, a.s., and therefore by standard stochastic integration:
 $$
 \dbE^0\Big[ \sum_{n=1}^{N} G_{t_{n-1}}(m_{t_{n-1}}, X_{t_{n-1}})
                                              \D^{\!\pi}\!X_{t_n}\Big] 
\underset{| \pi | \to 0}{\longrightarrow} \dbE^0\Big[\int_0^T G_{s-}(m_{s-}, X_{s-})dX_s \Big], \ \mbox{a.s.} $$
As for the remaining term, denoting $f:=\pa_x \d_m \phi,~g:=\pa_x \d_m \psi$, observe that:
 $$ 
 \big\{ G_{\cdot}(m_{t_{n-1}}, X_{t_{n-1}})\big\}_{t_{n-1}}^{t_n} 
 = \int_{t_{n-1}}^{t_n} f_{s}(m_{t_{n-1}}, X_{t_{n-1}})dA_s 
    + \int_{t_{n-1}}^{t_n} g_{s}(m_{t_{n-1}}, X_{t_{n-1}})dN_s.
$$
and let us dedicate the rest of this step to prove that: 
\begin{eqnarray}\label{CV1a}
\sum_{n=1}^N \D^{\!\pi}\!X_{t_n} \int_{t_{n-1}}^{t_n}f_{s}(m_{t_{n-1}},X_{t_{n-1}}) dA_s
&\underset{| \pi | \to 0}{\longrightarrow} &
\sum_{0 < s \le T} f_{s}(m_{s-}, X_{s-}) \D X_s \D A_s
\\
\mbox{and}~\sum_{n=1}^N \D^{\!\pi}\!X_{t_n} \int_{t_{n-1}}^{t_n}g_{s}(m_{t_{n-1}}, X_{t_{n-1}}) dN_s
&\hspace{-2mm}\underset{| \pi | \to 0}{\longrightarrow} &\hspace{-2mm}
\int_0^T g_{s}(m_{s-}, X_{s-})d[X, N]_s, \ \mbox{a.s.}~~~~
\label{CV1b}
\end{eqnarray}
To justify the convergence of \eqref{CV1a}, observe that:
\begin{eqnarray*}
\sum_{n=1}^{N} \D^{\!\pi}\!X_{t_n} 
                          \int_{t_{n-1}}^{t_n}f_{s}(m_{t_{n-1}},X_{t_{n-1}})
                                                       dA_s 
&=&  
\sum_{n=1}^{N} \D^{\!\pi}\!X_{t_n} \int_{t_{n-1}}^{t_n}\eta_{s, t_{n-1}} dA_s 
\\ 
&&+\sum_{n=1}^{N}  \D^{\!\pi}\!X_{t_n} \int_{t_{n-1}}^{t_n} f_{s}(m_{s-}, X_{s-})dA_s.
\end{eqnarray*}
where $\eta_{t_{n-1},s} := \{ f_{s}(m_{\cdot}, X_{\cdot})\}_{t_{n-1}}^s$. The first term in the right hand side of the last equation satisfies:
 \begin{eqnarray*}
\dbE\Big[ \Big| \sum_{n=1}^{N} \D^{\!\pi}\!X_{t_n}  \int_{t_{n-1}}^{t_n}  \eta_{s, t_{n-1}} dA_s \Big| \Big]  
&\le& 
\dbE\Big[\Big( \sum_{n=1}^{N} \big| \D^{\!\pi}\!X_{t_n} \big|^2 \Big)^{1/2} \Big(  \sum_{n=1}^{N} \big| \int_{t_{n-1}}^{t_n}  \eta_{s, t_{n-1}} dA_s \big|^2 \Big)^{1/2}\Big] 
\\
&\le& 
\dbE\Big[\sum_{n=1}^{N} \big| \D^{\!\pi}\!X_{t_n} \big|^2 \Big]^{1/2} \dbE\Big[\sum_{n=1}^{N} \big| \int_{t_{n-1}}^{t_n}  \eta_{s, t_{n-1}} dA_s \big|^2 \Big]^{1/2} 
\\
&\le& 
C\dbE\Big[\sum_{n=1}^{N} \big| \D^{\!\pi}\!X_{t_n} \big|^2 \Big]^{1/2} \dbE\Big[ \int_0^T  |\eta_{s, t_{n-1}(s)}|^2 d|A| _s \Big]^{1/2},
 \end{eqnarray*}
for some nonnegative constant $C$, where we used the boundedness of $A$. Note that $\dbE\big[\sum_{n=1}^{N} | \D^{\!\pi}\!X_{t_n} |^2 \big]$ is bounded due to the fact that $X$ has a finite quadratic variation. Moreover, as $\eta$ is bounded and $A$ has finite variation, it follows that $ \dbE\big[ | \int_0^T  \eta_{s, t_{n-1}(s)} dA_s |^2 \big] \longrightarrow 0$ as $| \pi | \to 0$, by the dominated convergence theorem for Stieltjes integrals. This shows that the last integral converges to $0$ as $| \pi | \to 0$ by dominated convergence. Therefore, by Lemma \ref{lem:cv-bracket} (i), we have in $\dbL^1$:
 \begin{align*}
  \lim_{| \pi | \to 0} \sum_{n=1}^{N} \D^{\!\pi}\!X_{t_n} \int_{t_{n-1}}^{t_n} f_{s}(m_{t_{n-1}}, X_{t_{n-1}})dA_s 
  &=  \lim_{| \pi | \to 0} \sum_{n=1}^N \D^{\!\pi}\!X_{t_n} \int_{t_{n-1}}^{t_n} f_{s}(m_{s-}, X_{s-})dA_s 
  \\ &= \int_0^T f_{s}(m_{s-}, X_{s-})d[X, A]_s  
  \\ &= \sum_{0 < s \le T} f_{s}(m_{s-}, X_{s-})\D X_s \D A_s,
 \end{align*}
 as $A$ is a finite variation process implying that $[X, A]^c = 0$. 
 
To justify the convergence of \eqref{CV1b}, observe that:
\begin{eqnarray*}
\dbE\Big[ \Big| \sum_{n=1}^N \D^{\!\pi}\!X_{t_n}  \int_{t_{n-1}}^{t_n} \tilde \eta_{t_{n-1}, s} dN_s \Big| \Big]  
&\le& \dbE\Big[\Big( \sum_{n=1}^N \big| \D^{\!\pi}\!X_{t_n} \big|^2 \Big)^{1/2} \Big(  \sum_{n=1}^N \big| \int_{t_{n-1}}^{t_n} \tilde \eta_{t_{n-1},s} dN_s \big|^2 \Big)^{1/2}\Big] 
\\
&\le&  \dbE\Big[\sum_{n=1}^N \big| \D^{\!\pi}\!X_{t_n} \big|^2 \Big]^{1/2} \dbE\Big[\sum_{n=1}^N \big| \int_{t_{n-1}}^{t_n} \tilde \eta_{t_{n-1},s} dN_s \big|^2 \Big]^{1/2}
\end{eqnarray*}
where $\tilde \eta_{t_{n-1},s} := \big\{ \pa_x \d_m \psi_{s}(m_{\cdot}, X_{\cdot}) \big\}_{t_{n-1}}^{s}$. Observe also that:
$$
\dbE\Big[ \sum_{n=1}^N \big| \int_{t_{n-1}}^{t_n} \tilde \eta_{t_{n-1},s} dN_s \big|^2 \Big] =  \dbE\Big[ \sum_{n=1}^N  \int_{t_{n-1}}^{t_n} | \tilde \eta_{t_{n-1},s} |^2 d[N]_s \Big] =  \dbE\Big[ \int_{0}^{T} | \tilde \eta_{t_{n-1}(s), s} |^2 d[N]_s \Big],
$$
and that the last term converges to $0$ as $| \pi | \to 0$, by dominated convergence. Then, up to some subsequence of partitions (still denoted the same), we have $ \dbE\Big[ \sum_{n=1}^N \big| \int_{t_{n-1}}^{t_n} \tilde \eta_{t_{n-1},s} dN_s \big|^2 \Big] \underset{| \pi | \to 0}{\longrightarrow} 0$, and therefore:
$$ \dbE\Big[ \Big| \sum_{n=1}^N \D^{\!\pi}\!X_{t_n}  \int_{t_{n-1}}^{t_n} \tilde \eta_{t_{n-1},s} dN_s \Big|\Big]  \underset{| \pi | \to 0}{\longrightarrow} 0.
$$
Since 
$$ \sum_{n=1}^N \D^{\!\pi}\!X_{t_n}  \int_{t_{n-1}}^{t_n} g_{s}(m_s, X_s) dN_s  \underset{| \pi | \to 0}{\longrightarrow} \int_0^T g_{s}(m_s, X_s)d[X, N]_s, \ \mbox{in $\dbL^1$}, $$
up to some subsequence, we may finally find a subsequence such that:
$$ \sum_{n=1}^N \D^{\!\pi}\!X_{t_n}  \int_{t_{n-1}}^{t_n} g_{s}(m_{t_{n-1}}, X_{t_{n-1}}) dN_s  \underset{| \pi | \to 0}{\longrightarrow} \int_0^T g_{s}(m_s, X_s)d[X, N]_s, \ \mbox{in $\dbL^1$,} $$
which implies the validity of \eqref{CV1b} along some further subsequence.
 \item For the second term in \eqref{terms}, we recall the notation from the previous step $G:=\pa_x \d_m U$ and we denote further $\hat F := \pa_{x\hat x}^2 \d_{mm}^2 U$. We first decompose:
 \begin{align*}
 \dbE^0\Big[ \{ G_{t_n}(\cdot, X_{t_{n-1},t_n}^{\l_2}) \}_{m_{t_{n-1}}}^{m_{t_{n-1}, t_n}^{\l_1}} \D^{\!\pi}\!X_{t_n}\Big] 
 =& \int_{[0,1]^2 }\dbE^0 \hat \dbE^0 \Big[\hat F_{t_n}(\theta_n^r)  \D^{\!\pi}\!X_{t_n} \D^{\!\pi}\!\hat X_{t_n} \Big]\l_1 d\l_4 d\l_5 \\
 =& \int_{[0,1]^2 }\dbE^0 \hat \dbE^0 \Big[ \hat F_{t_{n-1}}(\theta_n^0)\D^{\!\pi}\!X_{t_n} \D^{\!\pi}\!\hat X_{t_n} \Big]\l_4 d\l_4 d\l_5 \\
 &+ \int_{[0,1]^2 }\dbE^0 \hat \dbE^0 \Big[ \{\hat F\}^{t_n,\theta_n^r}_{t_{n-1},\theta_n^0} \D^{\!\pi}\!X_{t_n} \D^{\!\pi}\!\hat X_{t_n} \Big]\l_4 d\l_4 d\l_5, 
 \end{align*} 
 where we used the notations of Lemma \ref{lem:cv2} with $r=(\lambda_1\lambda_4,\lambda_2,\lambda_5)$.
 By Lemmas \ref{lem:cv-bracket} (i) and \ref{lem:cv2} (ii), we obtain, $\dbP$-almost surely:
 \bea\label{Lim1}
  &&\hspace{-5mm}
  \int_{[0,1]^2 }\dbE^0 \hat \dbE^0 \Big[ \sum_{n=1}^N \hat F_{t_n}(\th_n^r)  \D^{\!\pi}\!X_{t_n} \D^{\!\pi}\!\hat X_{t_n} \Big]\l_4 d\l_4 d\l_5 
  \\& &\hspace{-5mm}
  \underset{| \pi | \to 0}{\longrightarrow }\;
  \frac12 \dbE^0 \hat \dbE^0 \Big[ \int_0^T \hat F_{s-}(\th_s^0)d [X, \hat X]_s \Big] 
  \!+\!\! \int_{[0,1]^2} \!\! \dbE^0 \hat \dbE^0 \Big[\sum_{0 < s \le T} \{\hat F\}^{s,\theta^r_s}_{s-,\theta_{s}^0}  \D X_s \D \hat X_s  \Big] \l_4 d\l_4 d\l_5, 
  \nonumber
  \eea
\item  We now analyze the convergence of $\dbE^0\big[R_n\big]$ in \eqref{terms}. First, it follows from Lemma \ref{lem:cv2} (i)  that:
\begin{align*}
 \int_{[0,1]^2 }\dbE^0\Big[ \sum_{n=1}^{N} & \big\{F_{.}(m_{t_{n-1}},.)\big\}^{t_n,X_{t_{n-1}, t_n}^{\l_2 \l_3}}_{t_{n-1},X_{t_{n-1}}}| \D^{\!\pi} \!X_{t_n}|^2 \l_3 d\l_2 d\l_3\Big] \\
 &\underset{| \pi | \to 0}{\longrightarrow} 
 \!\int_{[0,1]^2} \!\!\dbE^0\Big[\sum_{0 < s \le T} \!\{F_.(m_{s-},.)\}^{s,X_{s}^{\l_2 \l_3}}_{s-,X_{s-}}|\D X_{s}|^2 \Big]\l_3 d\l_2 d\l_3 , \ \mbox{a.s.} 
\end{align*}
Now, observe that \eqref{eq:finiteness-derivatives} implies that $|F_{t_{n-1}}(m_{t_{n-1}}, X_{t_{n-1}})| \le C$ for some $C : \Omega \to \dbR_+$ measurable with respect to the sigma-field $\cF_T^0 \vee \sigma(A_{\cdot \wedge T}, N_{\cdot \wedge T})$. Finally, using Lemma \ref{lem:cv-bracket} for the conditional expectation $\dbE^0[\cdot]$, we deduce that
\begin{align*}
\sum_{n=1}^N \dbE^0\Big[ F_{t_{n-1}}(m_{t_{n-1}}, X_{t_{n-1}})| \D^{\!\pi} \!X_{t_n}|^2\Big] \underset{|\pi| \to 0}{\longrightarrow}  \dbE^0\Big[ \int_0^T F_{s-}(m_{s-}, X_{s-}) d[X]_s \Big], \q \mbox{a.s.}
\end{align*}
Finally, after rearranging all the terms with the right hand side of \eqref{Lim1} and separating discontinuities of the flow $\{m_t\}_{\{ t \in [0,T] \}}$, similarly to \eqref{jumpsXhatX} in the proof of Theorem \ref{thm:Ito}, we conclude that the desired formula holds true.\qed
\end{enumerate}

\end{theorem}

\section{Applications to stochastic control}\label{sec:examples}

\subsection{Mean field control with Poisson type common noise}

Let $\O := C^0([0,T], \dbR^d) \times \cP([0,T] \times \dbR_+ \times \dbR)$ be the canonical space, and $(X, \Pi^0)$ be the corresponding canonical process. We denote by $\dbF := \{\cF_t\}_{\{0 \le t \le T\}}$ the canonical filtration, and by $\dbF^0 := \{\cF_t^0\}_{\{0 \le t \le T\}}$ the filtration generated by $\Pi^0$, i.e., $\cF_t^0$ is the sigma-field generated by the restriction of $\Pi^0$ to $[0,t] \times \dbR_+ \times \dbR$. 
For a compact subset $\dbA \subset \dbR^k$, $k \in \dbN^*$, we denote by $\cA$ the set of $\dbF$-progressively measurable processes taking their values in $\dbA$. Given a measure $\dbP$ on $\O$, we denote $\dbP^0 := \dbP \circ (X | \Pi^0)^{-1}$. 

Let $b, \si, \g^0, \l^0$ and $f$ be bounded functions:
$$ (b, \si, \g, \l, f) : [0,T] \times \dbR^d \times \cP_2(\dbR^d) \times \dbA \longrightarrow \dbR^d \times \dbR^{d\times n} \times \dbR^{d_0} \times \dbR \times \dbR. $$
We assume that these functions are continuous in all their arguments, and Lipschitz-continuous in their $\dbR^d$ and $\cP_2(\dbR^d)$-valued arguments. In particular, we assume that the continuity in their third argument holds for both $\cW_1$ and $\cW_2$ distances. 
For $(t, \m) \in [0,T] \times \cP_2(\dbR^d)$, we denote by $\cP(t,\m)$ the set of $\dbP \in \cP_2(\O)$ such that: \\

${\rm{(i)}}$ $\Pi^0$ is a random Poisson measure with compensator ${\rm{Leb}}_{[0,T]} \otimes {\rm{Leb}}_{\dbR_+} \otimes \n(dy)$, \\ \\
${\rm{(ii)}}$ The process $X$ satisfies $\dbP \circ X_t^{-1} = \m$ and, for some $\a \in \cA$ and some $n$-dimensional $\dbP$-Brownian motion $W^\dbP$ independent of $\Pi^0$: 
$$ 
dX_s \!=\! b^{\a_s}(\Th_s)ds +  \si^{\a_s}(\Th_s)dW_s^\dbP +  \g^{\a_s}(\Th_s)dN_s, \ \,dN_s \!=\! \int_{_{\dbR_+ \!\times \dbR}} \!\!y \1_{_{\{\xi \le  \l^{\a_s}(\Th_s)\}}} \Pi^0(ds, d\xi, dy),
~\dbP-\mbox{a.s.} 
$$
for $s\ge t$, where $\Th_s := (s, \dbP^0_{X_s}, X_s)$ with $\dbP^0_{X_s}:=\dbP\circ(X_s|\cF^0_s)^{-1}$. 

\medskip
Our main interest is on the stochastic control problem:
\bea\label{MF-control-pb}
V(t, \m) := \sup_{\dbP \in \cP(t,\m)} \dbE^\dbP \Big[ \int_t^T f^{\a_s}\big(\Th_s \big)ds + g\big(\dbP^0_{X_T}\big)\Big],
\eea
where the control $\a \in \cA$ is associated with $\dbP$ accordingly to the definition above. 

Let ${\rm ST}^{\dbF^0}_{[t,T]}$ be the collection of all $\dbF^0-$stopping times with values in $[t,T]$. In this section, we assume that the following dynamic programming principle holds for all $t\in[0,T]$: 
\bea\label{MF-control-DPP}
&V(t, \m) = \sup_{\dbP \in \cP(t,\m)} \dbE^\dbP \Big[ \int_t^{\t^\dbP} f^{\a_s}\big(\Th_s \big)ds + V\big(\t^\dbP, \dbP_{X_{\t^\dbP}}^0\big)\Big], &
\\
&\mbox{for all family}~(\t^\dbP)_{\dbP\in\cP}\subset{\rm ST}^{\dbF^0}_{[t,T]},&
\nonumber\eea
see Djete, Possamaï \& Tan \cite{djete2022mckean} for specific conditions. Our goal is to derive the corresponding infinitesimal counterpart which reduces as standard to the so-called Hamilton-Jacobi-Bellman (HJB) equation. This requires to introduce the following integro-differential operators which characterize the infinitesimal generator of the controlled process and which are defined for all map $\f:\dbR_+\times\cP_2(\dbR^2)$ with appropriate smoothness by:
 \begin{align}
 &\cL^a \f := b^a \!\cdot\! \pa_x \d_m \f + \frac{1}{2}\si^a (\si^a)^\top\!:\!\pa_{xx}^2 \d_m \f, 
\label{La} \\
 &\cJ^a[\f](t,\m,x) :=  \l^a(t, \m, x) \int_{\dbR} \big\{ \f \big(t, \m \circ [\bx + y\g^a(t, \m, x)]^{-1}\big) - \f(t, \m)\big\} \n(dy),
\label{Ja}
\end{align}
  where $\bx$ is the identity map on $\dbR^d$.

\begin{proposition}\label{prop:HJB}
Let \eqref{MF-control-DPP} hold, and assume that $\pa_t V$, $\pa_x \d_m V$, $\pa_{xx}^2 \d_m V$, $\pa_{x \hat x}^2 \d_{mm}^2 V$ exist and are jointly continuous, with continuity in the measure argument in $\cW_1$ and $\cW_2$. Then $V$ is a classical solution of the HJB equation:
$$ 
-\pa_t V(t,\m) 
-  \int_{\dbR^d} \sup_{a \in \dbA} \big\{ f^a + \cL^a \d_m V 
                                                             +  \cJ^a[V] \big\}(t,\m,x) \m(dx)  = 0, 
~~(t, \m) \in [0,T] \times \cP_2(\dbR^d).
$$
\end{proposition}

\proof 
\textbf{Supersolution property.} Let $\bar \alpha_t = a \in \dbA$ be a constant control, and denote by $\bar \dbP$ the corresponding element of $\cP(t,\m)$. Set $\bar m_s := \bar \dbP_{X_s}^0$. By the dynamic programming principle \eqref{MF-control-DPP}, we have:
$$ V(t,\m) \ge \dbE^{\bar \dbP} \Big[ \int_t^{t+h} f^{\bar \alpha_s}\big(s, \bar m_s, X_s \big)ds + V\big(t+h, \bar m_{t+h}\big)\Big]. $$
By the same estimates as in Example \ref{ex:poisson}, we know that $X$ satisfies the estimates \eqref{integrability}, and therefore we obtain from Theorem \ref{thm:Ito} after appropriate localization:
\begin{equation}\label{supersolution}
\frac{1}{h}\int_t^{t+h} \!\!\dbE^{\bar \dbP}\Big[ - \pa_t V(s, \bar m_s) - (\cL^a \d_m V+ f^a)(s, \bar m_s, X_s)\Big]ds - \frac{1}{h} \dbE^{\bar \dbP}\Big[ \sum_{0 < s \le t} \D V(s, \bar m_s) \Big] \ge 0,~~
\end{equation}
where $\bar m_s := \bar \dbP^0_{X_{s}}$ and the operator $\cL^a$ is defined in \eqref{La}.
Now observe that:
\begin{eqnarray*}
\dbE^{\bar \dbP}\Big[ \sum_{t < s \le t+h} \D V(s, \bar m_s) \Big] 
&=& 
\dbE^{\bar \dbP}\Big[  \int_t^{t+h} \cJ^a[V](s, \bar m_s, X_s) ds \Big],
\end{eqnarray*}
where the integral operator $\cJ^a$ is defined in \eqref{Ja}. We plug this into \eqref{supersolution} and obtain the supersolution property by letting $h \to 0$, using the fact that $s \mapsto \bar m_s$ is right continuous and that $a \in \dbA$ is arbitrary. \\

\textbf{Subsolution property.}
Let $\tau^\dbP_h:=(t+h)\wedge\tau^\dbP$ with $\tau^\dbP:=\inf\{s>t:\cW_2(\m_s,\mu)\ge\delta\}$, for some $\d>0$, and $m_s := \dbP^0_{X_{s}}$.
By the dynamic programming principle \eqref{MF-control-DPP}, we have: 
$$
V(t, \m) \le \sup_{\dbP \in \cP(t,\m)} \dbE^\dbP \Big[ \int_t^{\tau^\dbP_h} f^{\a_s}\big(s, \dbP^0_{X_{s}}, X_s\big)ds + V\big(\tau^\dbP_h, \dbP_{X_{\tau^\dbP_h}}^0\big)\Big],
$$
from which we deduce, by It\^o's formula:
\begin{align*}
\inf_{\dbP \in \cP(t,\m)} \int_t^{\tau^\dbP_h} 
\dbE^{ \dbP}\Big[ &-\! (\pa_t V 
                               \!-\! \cL^{\a_s} \d_m V  
                               \!-\! f^{\a_s})(s,  \m_s,X_s) 
                               \!-\! \cJ^{\a_s}[V](s, m_{s}, X_{s-}) \Big]ds \le 0,
\end{align*}
For $\d > 0$, introduce now $B_2(\m,\d) := \{ \m' \in \cP_2(\dbR^d) : \cW_2(\m, \m') \le \d \}$ and $Q_2^{h,\d}(t,\m) := [t,t+h] \times B_2(\m,\d)$. Then it follows from the local boundedness of the derivatives of $V$ that: 
\begin{eqnarray*}
\inf_{(s,\m') \in Q_2^{h,\d}(t,\m)} \cH(s,\m')
\le
\frac{C_\delta}{h}\; \sup_{\dbP \in \cP(t,\m)} \int_t^{t+h} \dbP\big[m_s \notin B_2(\m, \d)\big]ds,
\end{eqnarray*}
for some constant $C_\d> 0$, where
$$
\cH(s,\m')
:=
- \pa_t V(s,  \m') -\int_{\dbR} \sup_{a \in \dbA} \Big\{ \cL^{a} \d_m V +f^a +\cJ^a[V]
                                                                                           \Big\}(s,\mu',x) \m'(dx).
$$
As the control process $\alpha$ take values in a the compact subset $\dbA$, it follows from the continuity of the coefficients of the controlled SDE that:
\begin{align*}
\dbP\big[ m_s \notin B_2(\m, \d)\big] =& \dbP\big[ \cW_2(m_s, \m) \ge \d \big] \le \frac{\dbE^\dbP\big[ | X_s - X_t |^2 \big]}{\d^2} \le C' \;\frac{s-t}{\d^2},
\end{align*}
for some constant $C'>0$. Together with the continuity of $\cH$, due to the compactness of $\dbA$, this implies that:
\begin{eqnarray*}
\inf_{\m' \in B_2(\m, \d)} \cH(t,\m')
\;=\;
\lim_{h\searrow 0}
\inf_{(s,\m') \in Q_2^{h,\d}(t,\m)}  \cH(s, \m')
&\le& 
\lim_{h\searrow 0}\frac{C'C_\d h}{2\d^2}
\;=\;
0.
\end{eqnarray*}
To conclude, we use the fact that, since all derivatives of $V$ involved in $\cH$ are in continuous both on $\cW_1$ and $\cW_2$, it follows that $\cH$ is uniformly continuous on $B_2(\m, 1)$ for the $\cW_1$ distance. Since $\cW_1 \le \cW_2$, we obtain by sending $\d\searrow 0$ that $\cH(t,\m)\le 0$, which is the required subsolution property. 
\qed

\subsection{Mean field optimal stopping with common noise}

We adapt the setting of Talbi, Touzi \& Zhang \cite{talbi2023dynamic} to the more general situation allowing for common noise. The canonical space is given by $\O := C^0([-1,T],\dbR^d) \times \dbI^0([-1,T]) \times C^0([-1, T], \dbR^{d_0})$, where $\dbI^0([-1,T])$ is the set of non-increasing and càdlàg maps from $[-1,T]$ to $\{0,1\}$,  constant  on $[-1,0)$, and ending with value $0$ at $T$. Note that the choice of the extension to $-1$ is arbitrary,  the extension of time  to the left of the origin is only needed to allow for an immediate stop at time $t=0$. We denote by $(X,I, W^0)$ the canonical process, with canonical filtration $\dbF = (\cF_t)_{t \in [-1,T]}$ and we denote by $\dbF^0$ the filtration generated by $W^0$. We also define $\mathbf{S}:=\dbR^d\times\{0,1\}$. 
\\

Let $b, \si, \g^0, \l^0$ and $f$ be bounded functions:
$$ (b, \si, \si^0, f) : [0,T] \times \dbR^d \times \cP(\dbR^d)  \longrightarrow \dbR^d \times \cM_{d, n} \times \cM_{d, d_0} \times \dbR. $$
As in the previous paragraph, we assume that these functions are continuous in all their arguments, and Lipschitz-continuous in their $\dbR^d$ and $\cP_2(\dbR^d)$-valued arguments.
For $(t, \m) \in [0,T] \times \cP_2(\bS)$, we denote by $\cP(t,\m)$ the set of $\dbP \in \cP_2(\O)$ such that: \\

\hspace{3mm} ${\rm{(i)}}$ $W^0$ is a $\dbP$-Brownian motion, 

\vspace{3mm}
\hspace{3mm} ${\rm{(ii)}}$ The process $X$ satisfies for some $n$-dimensional $\dbP$-Brownian motion $W^\dbP$ independent of $W^0$:
$$ 
\dbP \circ (X_t, I_{t-})^{-1} = m,
~\mbox{and}~
dX_s = I_s\Big(b(\Th_s)ds + \si(\Th_s)dW_s^\dbP + \si^0(\Th_s)dW_s^0\Big), ~s\ge t,~\dbP-\mbox{a.s.}
$$
where $\Th_s := (s, \dbP^0_{X_s}, X_s)$ and $\dbP^0_{X_s} := \dbP \circ (X_s | W^0)^{-1}$. 
\\

The value function of the mean field optimal stopping problem is defined by:
 \begin{equation}\label{weakoptstop}
 V(t,m) := \underset{\dbP \in \cP(t,m)}{\sup} \dbE^{\dbP}\Big[ \int_t^T f(\Th_r)dr  + g(\dbP_{X_T}^0) \Big], \quad \mbox{$(t,m) \in [0,T] \times \cP_2(\bS)$}.
 \end{equation}
Let ${\rm ST}^{\dbF^0}_{[t,T]}$ be the collection of all $\dbF^0-$stopping times with values in $[t,T]$. As in the previous paragraph, we assume that the following dynamic programming holds:
 \bea\label{MF-stopping-DPP}
 V(t,m) = \sup_{\dbP \in \cP(t,m)} \dbE^{\dbP}\Big[ \int_t^{\t^\dbP} f(\Th_r)dr  + V(\t, \dbP^0_{X_\t}) \Big],
 ~\mbox{for all family}~(\t^\dbP)\subset{\rm ST}^{\dbF^0}_{[t,T]}, 
 \eea
and we aim at deriving the dynamic programming equation as the infinitesimal counterpart of this dynamic programming principle.
\\

We recall from Talbi, Touzi \& Zhang \cite{talbi2023dynamic} the partial ordering on the Wasserstein space:
\begin{align*}
m' \preceq m
~\mbox{if, for some measurable}~p : \dbR^d \to [0,1],~&
m'(dx, 1) = p(x)m(dx, 1) ~\mbox{and}
\\& m'(dx, 0) = (1-p(x))m(dx, 1) + m(dx,0).
\end{align*}
Roughly speaking, $p$ corresponds to an immediate randomized stopping strategy sending the distribution $m$ to the distribution $m'$. We also introduce for all map $u:\dbR_+\times\cP_2(\dbR^d)\longrightarrow\dbR^d$ the notation
$$
C_u(t,m) := \{ m' \preceq m : u(t,m') = u(t,m) \}.
$$
Our objective is to prove that the dynamic programming equation is given by the obstacle problem on Wasserstein space:
\bea\label{obstacle}
\min_{m' \in C_u(t,m)} \!\big\{ \!-\!\mathbf{L} u (t, m') \big\} \!=\! 0, 
\,D_I u(t, m, \cdot)\!\ge\!0,
\,u(T, m) \!=\! g(m), \,(t,m) \!\in\! [0,T]\!\times\!\cP_2(\dbR^d),
\eea 
 where $D_I u(.,x) := \d_m u(.,x,1) - \d_m u(.,x,0)$ and, denoting $\cL^0:=b\!\cdot\!\pa_x+ \frac{1}{2}\si \si^\top\!:\!\pa_{xx}^2$, the operator $\mathbf{L}$ is defined by
 \begin{eqnarray*}
 \mathbf{L} u(t,m)
 &\hspace{-2mm}:=& \hspace{-2mm}
 \pa_t u(t,m) 
 \\
 &&\hspace{-10mm}
 +\!\int_{\dbR^d}\!\! \Big[ f\!+\!\cL^0 \d_m u(.,\!1)
                                              \!+\!\frac{1}{2}\int_{\dbR^d}\!\! \si^0 {\si^0}^\top\!\!\!:\!
                                                                                     \pa_{x \hat x}^2 \d_{mm}^2 u(.,1, \hat x,\!1)) m(d\hat x,\!1)
                                     \Big](t,m,x)m(dx,\!1).
 \end{eqnarray*}
We refer to \cite{talbi2023dynamic, talbi2023viscosity} for more comments about the structure of the obstacle problem \eqref{obstacle} in the context of mean field optimal stopping without common noise. 
 
 \begin{proposition}\label{prop:obstacle}
 Let \eqref{MF-stopping-DPP} hold, and assume that $\pa_t V$, $\d_m V$, $\pa_x \d_m V$, $\pa_{xx}^2 \d_m V$ and $\pa_{x \hat x}^2 \d_{mm}^2 V$ exist and are continuous, with bounded $\pa_{xx}^2 \d_m V$ and $\pa_{x \hat x}^2 \d_{mm}^2 V$. Then $V$ is a classical solution of the obstacle problem \eqref{obstacle}.
\end{proposition}
\proof
By taking $\t =t$ and $\dbP$ such that $\dbP_{(X_t, I_t)} = \dbP^0_{(X_t, I_t)} = m' \preceq m$ in the dynamic programming principle \eqref{MF-stopping-DPP}, we obtain: $V(t,m) \ge V(t,m')$. Since this is true for any $m' \preceq m$, this implies that $D_I V(t,m,\cdot) \ge 0$ for the same reasons as in \cite[Lemma 4.3]{talbi2023dynamic}. \\

Now fix $h > 0$ and define $\th_h^\dbP := \inf \{ s \ge t : \cW_2\big(\dbP_{(X_t, I_t)}, \dbP_{(X_s, I_s)}\big) \ge 1\} \wedge (t+h)$. Let $\bar \dbP \in \cP(t,m)$ such that $I$ is constant under $\bar \dbP$. By \eqref{MF-stopping-DPP} and Theorem \ref{thm:Ito}, we have:
\begin{equation}\label{obstacle-ineq1}
-\frac{1}{h}\int_t^{\th_h^{\bar \dbP}} \dbL u(s, \bar m_s)ds \ge 0,
\end{equation}
where $\bar m_s := \bar \dbP^0_{(X_s, I_s)}$. Here, we used the fact that the stochastic integral w.r.t. \ $W^0$ vanishes due to the boundedness of $\si \si^\top \pa_x \d_m V(s, \bar m_s, X_s)$ is bounded on $[t, \th_h]$, $\bar \dbP$-a.s., and the process $I$ is constant, $\bar \dbP$-a.s. implying that the contributions of the jumps equal to $0$. Now, by continuity of the flow of marginals $s \mapsto \bar m_s$, it is clear that $\th_h = t+h$ for $h$ smaller than some (random) $\bar h$. Therefore, sending $h \to 0$ in \eqref{obstacle-ineq1}, we obtain:
$$ -\dbL V(t,m) \ge 0. $$ 
Now, observe that, as in \cite[Proposition 2.2]{talbi2023dynamic}, the set $\cP(t,m)$ is compact for the 2-Wasserstein distance. Thus, by continuity of $F$ and $V$, there exists $\hat \dbP \in \cP(t,m)$ such that:
$$ V(t,m) = \dbE^{\hat \dbP}\Big[ \int_t^{\th_h^{\hat \dbP}} F(s, \hat m_s)ds + V\big(\th_h^{\hat \dbP}, \hat m_{\th_h^{\hat \dbP}}\big) \Big], $$
where $\hat m_s := \hat \dbP^0_{(X_s, I_s)}$. Then, by Theorem \ref{thm:Ito} again, we have:
\begin{align*}
- \int_t^{\th_h^{\hat \dbP}} \dbL V(s, \hat m_s)ds - \dbE^{\hat \dbP}\Big[\int_t^{\th_h^{\hat \dbP}} D_I V(s, \hat m_s, X_s)dI_s + \sum_{t \le s \le \th_h^{\hat \dbP}} \D V(s,\hat m_s) \Big] = 0.
\end{align*}
We have proved that those three terms are nonnegative, hence:
$$ -\frac{1}{h} \int_t^{\th_h^{\hat \dbP}} \dbL V(s, \hat m_s)ds = 0. $$
Since $s \mapsto \hat m_s$ is right continuous, $\hat \dbP$-a.s., we obtain by sending $h \to 0$:
$$ -\dbL V(t, \hat m_t)= 0. $$
We conclude by observing that $\hat m_t \in C_V(t,m)$. 
\qed

\appendix

\section{Technical Results}\label{sec:appendix}

\subsection{Proof of Lemma \ref{lem:cv-bracket}}
(i) As in the previous proofs, we show the result in the one dimensional setting $d =1$. Observe that:
$$  \sum_{n=1}^N H_{t_{n-1}}^\pi | \D^{\!\pi}\!X_{t_n} |^2 
=  \sum_{n=1}^N H_{t_{n-1}}^\pi \D^{\!\pi} \!\big(X^2 \big)_{t_n} 
    - 2 \sum_{n=1}^N H_{t_{n-1}}^\pi X_{t_{n-1}} \D^{\!\pi}\!X_{t_n}. $$
We are then reduced to prove the $\dbL^1$ convergence of integrals with respect to a semimartingale. Let us detail the argument for the first term in the right hand side equality. We prove the convergence of the finite variation part and the local martingale part independently. First observe that:
$$ 
X^2 = X_0^2 + \tilde V + \tilde M
~\mbox{with}~
\tilde V_t := 2\int_0^t X_s dV_s + [M]_t, \q \tilde M_t := 2\int_0^t X_s dM_s,~t\ge 0. 
$$
Note that:
$$ \Big| \sum_{n=1}^N H_{t_{n-1}}^\pi  \D^\pi \tilde V_{t_n} - \int_0^T H_{s-}d\tilde V_s \Big| \le \int_0^T \big| H_{t_{n-1}(s)}^\pi - H_{s-} \big| d| \tilde V|_s. $$
As $H_{t_{n-1}(s)}^\pi - H_{s-} \longrightarrow 0$ a.s., the $\dbL^1-$convergence to zero of the last integral follows from the dominated convergence theorem due to the boundedness of $H$ together with the observation that $\sup_{ t \in [0,T]} | X_t |^2  \in \dbL^1$ by \eqref{integrability} and the BDG inequality, implying that:
$$ 
| \tilde V |_T \le 2 \sup_{ t \in [0,T]} | X_t | | V |_T + [M]_T \le  \sup_{ t \in [0,T]} | X_t |^2 +  | V |_T^2 + [M]_T \in \dbL^1.
$$

As for the martingale part, using BDG inequality, we have:
\begin{eqnarray*} 
\dbE\Big[\Big| \sum_{n=1}^N H_{t_{n-1}}^\pi \D^\pi\tilde M_{t_n} - \int_0^T H_{s-} d\tilde M_s \Big|\Big] 
&\le& \dbE\Big[ \Big( \int_0^T \big| H_{t_{n-1}(s)}^\pi  - H_{s-}  \big|^2 d[ \tilde M ]_s \Big)^{1/2}\Big] 
\\
&=& 
2 \dbE\Big[ \Big( \int_0^T \big| H_{t_{n-1}(s)}^\pi  - H_{s-}  \big|^2 X_s^2 d[ M ]_s \Big)^{1/2}\Big]. 
\end{eqnarray*}
Again, since $H^\pi - H$ is uniformly bounded, there exists a constant $C \ge 0$ such that:
$$ \Big( \int_0^T \big| H_{t_{n-1}(s)}^\pi  - H_{s-}  \big|^2 X_s^2 d[ M ]_s \Big)^{1/2} \le C \Big(\sup_{t \in [0,T]} X_t^2 [M]_T\Big)^{1/2} \in \dbL^1 $$
by Cauchy-Schwarz inequality. We may then apply again the dominated convergence theorem to derive the desired result. 

\vspace{3mm}
\noindent (ii) 
Denote $R_n^\pi := \{ f \}_{\th_n^0}^{\th_n^r} \D^{\!\pi} X_{t_n} \D^{\!\pi} \hat X_{t_n}$, $R_s := \{ f \}_{\th_s^0}^{\th_s^r}\D X_s \D \hat X_s$, and:
$$ \bar R_\pi := \sum_{n=1}^{N}R_n^\pi, \ \bar R := \sum_{0 < s \le T} R_s. $$
By continuity of $f$ in its arguments, $f$ is uniformly bounded on any bounded set of $\cP_2(\dbR^d) \times \dbR^d \times \dbR^d$ (where the first coordinate is measured with the 1-Wasserstein distance). Therefore, there exists a modulus of continuity $\rho$ such that:
\bea\label{modulus}
 | R_n^\pi | \1_{\{ | \xi_{n-1, n} | \le 1\}} \le \rho\big(\xi_n^\pi \big), \q \mbox{with} \q \xi_n^\pi :=  \cW_1(m_{t_{n-1}}, m_{t_{n}}) + | \D^{\!\pi}\!X_{t_n} | + | \D^{\!\pi}\!\hat X_{t_n}|
 \eea
Now given $\eta \in (0,1)$, we consider a continuous function $h_\eta : \dbR_+ \to [0,1]$ such that $h_\eta = 1$ on $[2\eta, \infty)$ and $h_\eta = 0$ on $[0, \eta]$. First observe that, for $| \pi |$ sufficiently small, the number of terms such that $| R_n^\pi |  h_\eta(\xi_n^\pi) \neq 0$ is finite, as $\sum_{n = 1}^N | \D^{\!\pi} X_{t_n} | \longrightarrow \sum_{s \le T} | \D X_s | < \infty$ a.s.\ as $\pi \to 0$ by the summability conditions \eqref{integrability}. Then it follows from the continuity of $h_\eta$ that
\begin{equation}\label{finite-sum} 
 \sum_{n=1}^{N} R_n^\pi h_\eta(\xi_n^\pi) \\ \underset{| \pi | \to 0}{\longrightarrow}  \sum_{0 < s \le T} R_s h_\eta(\xi_s), \q \mbox{where} \q \xi_s := \cW_1(m_{s-}, m_s) + | \D X_s | + | \D \hat X_s |.
\end{equation}
as these two sums are in fact finite a.s. As for the small jumps:
\begin{align*}
\Big\lvert  \sum_{n=1}^{N} R_n^\pi[1-h_\eta(\xi_n^\pi)] \Big\rvert &\le  \rho(\eta)  \sum_{n=1}^{N} \big\lvert \D^{\!\pi}\!X_{t_n} \D^{\!\pi}\!\hat X_{t_n} \big\rvert \\
&\le  \frac12\rho(\eta)({\rm QV}_{\!\pi}(X) + {\rm QV}_{\!\pi}(\hat X)), 
\end{align*}
where, for any process $Y$, we write ${\rm QV}_{\!\pi}(Y) := \sum_{n=1}^{N} | \D^\pi Y_{t_{n}} |^2$. 
As the sum on the right hand side converges in probability to the quadratic covariation of $X$ and $\hat X$, and $\eta$ may be chosen arbitrarily small, we deduce that the left hand side converges to $0$ in probability as $| \pi | \to 0$. Finally, combining this result with \eqref{finite-sum}, we obtain:
\begin{eqnarray*} 
| \bar R_\pi \!-\! \bar R| 
&\hspace{-2mm}\le& \hspace{-2mm}
\big\lvert  \sum_{n=1}^{N} R_n^\pi h_\eta(\xi_n^\pi) 
- \!\! \sum_{0 < s \le T}\!\! R_s h_\eta(\xi_s) \big\rvert  + \frac12 \rho(\eta) ({\rm QV}_{\!\pi}(X) \!+\! {\rm QV}_{\!\pi}(\hat X)) 
     +  \!\sum_{0 < s \le T} \lvert R_s \rvert \1_{\{\xi_s \le \eta \}} \\
\\
&&
\underset{| \pi | \to 0}{\longrightarrow} \frac12 \rho(\eta) \big( [X]_T + [\hat X]_T \big) +   \sum_{0 < s \le T} \lvert R_s \rvert \1_{\{\xi_s \le \eta \}} \le \rho(\eta) \big( [X]_T + [\hat X]_T \big),
\end{eqnarray*}
where the convergence holds in probability. Sending $\eta$ to $0$, this shows that $| \bar R_\pi - \bar R|$ converges to $0$ in probability as $| \pi | \to 0$. Since  $\bar R_\pi \le \frac12 | f |_\infty \big({\rm QV}_{\!\pi}(X) \!+\! {\rm QV}_{\!\pi}(\hat X)\big)$, which converges in $\dbL^1$ by Lemma \ref{lem:cv-bracket} (i), we see that the family $\{| \bar R_\pi - r |\}_\pi$ is uniformly integrable, and therefore we deduce that $\bar R_\pi \longrightarrow \bar R$ in $\dbL^1$. 
\qed

\subsection{Proof of Lemma \ref{lem:cv2}}
{\rm (i)} We first decompose:
\begin{align*}
 &\sum_{n=1}^{N} \big\{F(m_{t_{n-1}},.)\big\}^{t_n,X_{t_{n-1}, t_{n}}^\l}
                                                                      _{t_{n-1},X_{t_{n-1}}}
                             | \D^{\!\pi}\!X_{t_{n}}|^2 = \bar R_\pi + \bar S_\pi,
\end{align*}
where:
\begin{align*}
\bar R_\pi := \sum_{n=1}^{N} \big\{ F_{t_{n}}(m_{t_{n-1}}, \cdot) \big\}_{X_{t_{n-1}}}^{X_{t_{n-1}, t_{n}}^\l} | \D^{\!\pi}\!X_{t_{n}}|^2, \q \bar S_\pi := \sum_{n=1}^{N} \big\{ F_{\cdot}(m_{t_{n-1}}, X_{t_{n-1}}) \big\}_{t_{n-1}}^{t_{n}} | \D^{\!\pi}\!X_{t_{n}}|^2.
\end{align*}
We shall prove that the following convergence results hold:
\begin{eqnarray}
\bar R_\pi 
\underset{| \pi | \to 0}{\longrightarrow} 
\bar R:= \sum_{0 < s \le t} \big\{ F_{s}(m_{s-}, \cdot) \big\}_{X_{s-}}^{X_{s}^\l} | \D X_{s} |^2,~~
\mbox{in}~\dbL^1(\dbP_{| \cF_T^0, A, N}),
\label{R_pi}
\\
\bar S_\pi 
\underset{| \pi | \to 0}{\longrightarrow} 
\bar S:= \sum_{0 < s \le t} \big\{ F_{\cdot}(m_{s-}, X_{s-}) \big\}_{s-}^{s} | \D X_{s} |^2
,~~
\mbox{in}~\dbL^1(\dbP_{| \cF_T^0, A, N}),
\label{S_pi}
\end{eqnarray} 
where $\dbP_{| \cF_T^0, A, N}$ denotes the conditional distribution of the canonical process given $\cF_T^0, A$ and $N$. 
\begin{itemize}
\item \textit{Convergence of $\bar R_\pi$}. For $\eta \in (0,1)$, using the function $h_\eta$ introduced in the proof of Lemma \ref{lem:cv-bracket}, we have for some (random) modulus of continuity $\rho$: 
 \begin{align*} 
\sum_{n=1}^{N} | \D^{\!\pi}\!X_{t_{n}} |^2  \big\{ F_{t_{n-1}}(m_{t_{n-1}}, \cdot) \big\}_{X_{t_{n-1}}}^{X_{t_{n-1}, t_{n}}^\l}  [1- h_\eta(| \D^{\!\pi}\!X_{t_{n}} |)]  \le \rho(\eta) \sum_{n=1}^N | \D^{\!\pi}\!X_{t_{n}} |^2.
\end{align*} 
  Similarly, we have:
$$ \sum_{0 < s \le T} | \D X_{s} |^2  \big\{ F_{s-}(m_{s-}, \cdot) \big\}_{X_{s-}}^{X_{s}^\l} [1- h_\eta(| \D X_{s} |)] \le \rho(\eta) \sum_{0 <s \le T} | \D X_{s} |^2. $$ 
By the summability conditions on $X$, the number of $s \in [0,T]$ such that $| \D X_s | \ge \eta$ is finite, and therefore, as both sums are finite,
\begin{align*}
\Big| \sum_{n=1}^{N} | \D^{\!\pi}\!X_{t_{n}} |^2 
&\big\{ F_{t_{n-1}}(m_{t_{n-1}}, \cdot) \big\}_{X_{t_{n-1}}}^{X_{t_{n-1}, t_{n}}^\l}  
h_\eta(| \D^{\!\pi}\!X_{t_{n}} |) 
\\
& - \sum_{0 < s \le T} | \D X_{s} |^2 \big\{ F_{s-}(m_{s-}, \cdot) \big\}_{X_{s-}}^{X_{s}^\l} 
h_\eta(| \D X_s |)  \Big| 
\;\underset{| \pi | \to 0}{\longrightarrow}\; 0, \ \dbP-\mbox{a.s.} 
\end{align*}
 By the same reasoning as in the proof of Lemma \ref{lem:cv-bracket} (ii), we deduce that $\bar R_\pi$ converges to $\bar R$ in probability. Moreover, \eqref{eq:finiteness-derivatives} implies that $F$ is bounded by a random variable $C$ measurable with respect to $\cF_T^0 \vee \sigma(A_{\cdot \wedge T}, N_{\cdot \wedge T})$. Since $\mathrm{QV}_\pi(X)$ and $\mathrm{QV}(X)$ are uniformly integrable,  we deduce that $\bar R_\pi$ and $\bar R$ are uniformly integrable for the conditional expectation $\dbE^0$, and therefore that \eqref{R_pi} holds.
\item \textit{Convergence of $\bar S_\pi$}. Let $F^c$ and $F^a$ respectively denote the continuous and atomic parts of $F$. By pathwhise continuity of $t \mapsto F_t^c$ on the compact $[0,T]$, there exists a (random) modulus of continuity such that:
\begin{align*} 
\Big| \sum_{n=1}^N | \D^{\!\pi}\!X_{t_n} |^2 &\big\{ F_{\cdot}^c(m_{t_{n-1}}, X_{t_{n-1}}) \big\}_{t_{n-1}}^{t_n} 
\\ &
- \sum_{0 < s \le T} | \D \!X_s |^2 \big\{ F_{\cdot}^c(m_{s-}, X_{s-}) \big\}_{s-}^s \Big| \le \rho(| \pi |)\Big( \mathrm{QV}_\pi(X) + \mathrm{QV}(X)\Big).    
\end{align*}   
Furthermore, we have:
\begin{align*} 
\Big| \sum_{n=1}^N | \D^{\!\pi}\!X_{t_n} |^2 \big\{ F_{\cdot}^a(m_{t_{n-1}}, X_{t_{n-1}}) \big\}_{t_{n-1}}^{t_n} - \sum_{0 < s \le T} | \D \!X_s |^2 \big\{ F_{\cdot}^a(m_{s-}, X_{s-}) \big\}_{s-}^s \Big| \underset{| \pi | \to 0}{\longrightarrow} 0, \q \mbox{a.s.}    
\end{align*}   
Therefore, we have $\bar S_\pi \to \bar S$, a.s., as $| \pi | \to 0$. As in the case of $\bar R_\pi$, the sequence $\bar S_\pi$ is uniformly integrable for the conditional expectation $\dbE^0$, which implies that
\eqref{S_pi} holds.  
\end{itemize}
{\rm (ii)}  Observe that:
\begin{align*}
 &\sum_{n=1}^{N} \big\{\hat F_.(\cdot))\big\}^{t_n,\th_n^\l}
                                                                      _{t_{n-1},\th_n^{\bm{0}}}
                             | \D^{\!\pi}\!X_{t_{n}}|^2 = \bar R_\pi + \bar S_\pi,
\end{align*}
where:
\begin{align*}
\bar R_\pi := \sum_{n=1}^{N} \big\{ \hat F_{t_{n}}(\cdot) \big\}_{\th_n^{\bm{0}}}^{\th_n^\l} | \D^{\!\pi}\!X_{t_{n}}|^2, \q \bar S_\pi := \sum_{n=1}^{N} \big\{ \hat F_{\cdot}(\th_n^{\bm{0}}) \big\}_{t_{n-1}}^{t_{n}} | \D^{\!\pi}\!X_{t_{n}}|^2.
\end{align*} 
We need to prove that:
\begin{eqnarray}
\bar R_\pi  
\underset{| \pi | \to 0}{\longrightarrow} 
\bar R:= \sum_{0 < s \le t} \big\{ \hat F_{s}( \cdot) \big\}_{\th_s^{\bm{0}}}^{\th_s^{\l}} | \D X_{s} |^2,~~
\mbox{in}~\dbL^1(\dbP_{| \cF_T^0, A, N}),
\label{2R_pi}
\\
\bar S_\pi 
\underset{| \pi | \to 0}{\longrightarrow} 
\bar S:= \sum_{0 < s \le t} \big\{ \hat F_{\cdot}(\th_s^{\bm{0}}) \big\}_{s-}^{s} | \D X_{s} |^2
,~~
\mbox{in}~\dbL^1(\dbP_{| \cF_T^0, A, N}).
\label{2S_pi}
\end{eqnarray}
\eqref{2R_pi} follows from the same arguments as the convergence of $\bar R_\pi$ in (i), replacing $| \D^\pi X_{t_n} |$ with
$$ \xi_n^\pi :=  \cW_1(m_{t_{n-1}}, m_{t_{n}}) + | \D^{\!\pi}\!X_{t_n} | + | \D^{\!\pi}\!\hat X_{t_n}|$$
as argument of the function $h_\eta$, and using the fact that $\hat f$ and $\hat g$ are continuous with respect to both distances $\cW_2$ and $\cW_1$ to guarantee the existence of a modulus of continuity $\rho$. \eqref{2S_pi} follows from the very same arguments as (i). 
\qed 

\bibliography{Bibliography}

\end{document}